\documentclass[12pt]{article}
\RequirePackage{amsthm,amsmath}
\usepackage{amsfonts}

\usepackage{color,latexsym,amsfonts,amssymb}

\newcommand{\ignore}[1]{}{}
\newcommand{\be}{\begin{equation}}               
\newcommand{\ee}{\end{equation}}                 
\newcommand{\bi}{\begin{itemize}}
\newcommand{\ei}{\end{itemize}}

\newcommand{\beqn}{\begin{eqnarray}}             
\newcommand{\eeqn}{\end{eqnarray}}               
\newcommand{\beq}{\begin{eqnarray*}}             
\newcommand{\eeq}{\end{eqnarray*}}               

\newcommand{\nn}{\nonumber}

\newcommand{\lbl}{\label}

\newcommand{\ssb}{\scriptstyle \footnotesize 
                 \begin{array}{c}}

\newcommand{\esb}{\end{array}}

\begin{document}
\baselineskip=7.0mm

\renewcommand{\theequation}{\arabic{section}.\arabic{equation}}
\newtheorem{theorem}{{\sc Theorem}}[section]
\newtheorem{prop}{ {\sc Proposition}}[section]
\newtheorem{lemma}{{\sc Lemma}}[section]
\newtheorem{coro}{{\sc Corollary}}[section]
\newtheorem{remark}{{\sc Remark}}[section]
\newtheorem{exam}{{\sc Example}}[section]

\newcommand{\Proof}{ \noindent {\bf Proof.}\ }

\begin{flushleft}
\hspace{-.3cm} {\bf \large On the longest length of  arithmetic progressions }

\bigskip

Min-Zhi Zhao\footnote{Research supported  by NSFC (Grant No. 11101113 ) and  by ZJNSF (Grant No. R6090034)}, Department of Mathematics, Zhejiang University, Hangzhou, Zhejiang, China\\
{\it E-mail}: zhaomz@zju.edu.cn

\medskip
Hui-Zeng Zhang \footnote{Research supported  by NSFC (Grant No. 11001070)},
Department of Mathematics,  Hangzhou Normal University, Hangzhou, Zhejiang, China\\
{\it E-mail}: zhanghz789@163.com

\end{flushleft}

\vspace{1cm}

{\small
{\bf Abstract.}
Suppose that $\xi^{(n)}_1,\xi^{(n)}_2,\cdots,\xi^{(n)}_n$ are i.i.d with $P(\xi^{(n)}_i=1)=p_n=1-P(\xi^{(n)}_i=0)$.
Let $U^{(n)}$  and $W^{(n)}$  be the   longest length of arithmetic progressions and of
arithmetic progressions  mod $n$ relative to  $\xi^{(n)}_1,\xi^{(n)}_2,\cdots, \xi^{(n)}_n$
respectively. Firstly, the  asymptotic distributions  of $U^{(n)}$ and $W^{(n)}$ are given.  Simultaneously,
the errors  are estimated by using Chen-Stein method.
Next, the almost surely limits are discussed when all $p_n$ are equal and when considered on a common probability space.
Finally, we consider the case that $\lim_{n\to\infty}p_n=0$ and
$\lim_{n\to\infty}{np_n}=\infty$.  We prove that as $n$ tends to $\infty$,
 the probability that
 $U^{(n)}$ takes two numbers  and $W^{(n)}$  takes three numbers tends to $1$ .
}

\vspace{1cm}


\vspace{.6cm}


\newpage

\section{Introduction and main results} \lbl{sect1}

 Limit distributions for  the longest length of  runs with respect to Bernoulli sequence   have
been investigated for a long time, see, e.g.,  \cite{ER70}, \cite{ER77},\cite{GSW86}, \cite{J09} and \cite{M93}.
But what about the longest length of arithmetic  progressions?
Problems   connected to arithmetic  progressions
 are very important in number theory, see \cite{T07}. For example, Roth's theory says that
every set of integers of positive density contained infinitely many progressions of length three.

Suppose that $\xi_1,\xi_2,\cdots$ is a Bernoulli sequence  with $P(\xi_i=1)=p=1-q$, where
$0<p<1$.
Let $\Sigma_n=\{1\le i\le n:\xi_i=1\}$ be the random subset of $\{1,2,\cdots n\}$ decided by
$\xi_1,\xi_2,\cdots,\xi_n$.
 For any $1\le a,s\le n$ ,
define
$$U^{(n)}_{a,s}=\max\{1\le m\le 1+[\frac {n-a}{s}]:\xi_a=1,\xi_{a+s}=1,\cdots,\xi_{a+(m-1)s}=1 \},$$
which is the maximum length of  arithmetic progressions in $\Sigma_n$ starting at $a$, with difference $s$.
Let
$$U^{(n)}=\max_{1\le a,s\le  n} U^{(n)}_{a,s},$$
which is the length of the longest  arithmetic progression in $\Sigma_n$. We call $U^{(n)}$ the longest  length
of  arithmetic progressions relative to $\xi_1,\xi_2,\cdots,\xi_n$.

For any $1\le a,s\le n$, the numbers
$$a, a+s \mod n, a+2s \mod n, \cdots, a+({n}/{\gcd(s,n)}-1) s\mod n$$ are different while
$a+\big( {n}/{\gcd(s,n)}\big ) s \mod n=a,$ where $\gcd(s,n)$ denotes the greatest  common divisor of $s$ and $n$.
For convenience, let $kn \mod n=n$ for any integer $k$.
Define
$$W^{(n)}_{a,s}=\max\{1\le m\le \frac {n}{\gcd (s,n)} :\xi_a=1,\xi_{a+s \mod n}=1,\cdots,\xi_{a+(m-1)s\mod n}=1 \},$$
and
$$W^{(n)}=\max_{1\le a,s\le n} W^{(n)}_{a,s}.$$
We call $W^{(n)}$ the  the longest length
of  arithmetic progressions mod $n$ relative to $\xi_1,\xi_2,\cdots,\xi_n$.
Note that  $U^{(n)}$ is an increasing  of $n$ while
$W^{(n)}$ is not.

In \cite{BYZ07}, the authors discussed the limit distribution
of $U^{(n)}$ and $W^{(n)}$ in the case that $p=1/2$.
The results can be easily  extended to the case that $p\not=1/2$.
Set $C= {-2}/{\ln p}$ and let $\ln$ denote the logarithm of base $e$.  In  \cite{BYZ07}, they proved that
as $n$ tends to $\infty$,
$\frac {U^{(n)}}{C\ln n}\rightarrow 1 $  in probability  and
$\frac {W^{(n)}}{C\ln n}\rightarrow 1 $ in probability. Furthermore,
 $\lim_{n\to\infty} \frac {U^{(n)}}{C\ln n}=1 \,\, a.s.$ and
\beqn
1=\liminf_{n\to\infty}\frac {W^{(n)}}{C\ln n}<\frac 32\le \limsup_{n\to\infty}\frac {W^{(n)}}{C\ln n}\,\, a.s.\lbl{d1}
\eeqn
The authors also conjectured that
\beqn
 \limsup_{n\to\infty}\frac {W^{(n)}}{C\ln n}=\frac 32 \,\,a.s.\lbl{d2}
\eeqn

In this paper, we will use Chen-Stein method to study the asymptotic distributions
of $U^{(n)}$ and $W^{(n)}$ more carefully. In addition, the errors  are also given.  The limit distributions we get
are a bit  different  from that in \cite{BYZ07}.
Next, we {\bf prove the conjecture } and give more description
about  almost surely limits.

Set  $D=1/{\ln p}$.
 For any integer $1<r\le n$, let
\begin{eqnarray}
\lambda_{n,r}=\lambda_{n,r,p}=\frac {n^2p^r(p+qr)}{2r(r-1)}.\lbl{s1-d1}
\end{eqnarray}
and
\begin{eqnarray}
 \mu_{n,r}=\mu_{n,r,p}=\frac {qn^2p^r}{2}. \lbl{s1-d2}
\end{eqnarray}

\begin{theorem} \lbl{t1}
{\rm (1)} It holds that
\begin{eqnarray}
 \max_{1< r\le n} |P(U^{(n)}<r)-e^{-\lambda_{n,r}}|\le O(\frac {\ln^{4} n\ln\ln n}{n}).\lbl{th1-d1}
\end{eqnarray}

{\rm (2)} Let $h_n=C\ln n+ D\ln\ln n$.  For any $a<1$,
\begin{equation}
 \max_{x\ge  aD\ln\ln n,h_n+x\in \mathbb Z}|\exp(\frac {-q\ln p}{4} p^{x}) P(U^{(n)}<h_n +x)-1|\le O(\frac{\ln\ln n}{\ln^{1-a} n}).\lbl{th1-d3}
\end{equation}

{\rm(3)} As $n$ tends to $\infty$,
\begin{eqnarray}
\frac {U^{(n)}-C\ln n}{\ln \ln n}\rightarrow D \text{ in probability }.\lbl{th1-d2}
\end{eqnarray}

{\rm(4)}Almost surely,
\beqn
&&D=\liminf_{n\to\infty}\frac {U^{(n)}-C\ln n}{\ln\ln n}<\limsup_{n\to\infty}\frac {U^{(n)}-C\ln n}{\ln\ln n}=0,\lbl{th1-d4}\\
&&\lim_{n\to\infty}\frac 1{n\ln n}(\sum_{k=1}^n U^{(2^k)}+ \frac{\ln 2}{\ln p} n^2)=D\lbl{th1-d5}
\eeqn
and
\beqn
\lim_{n\to\infty}\frac 1 {\ln n\ln\ln n}(\sum_{k=1}^n \frac {U^{(k)}}{k}+ D\ln^2 n)=D.\lbl{th1-d6}
\eeqn
\end{theorem}

\begin{theorem} \lbl{t2}
{\rm(1)}It holds that
\begin{eqnarray}
 \max_{1< r\le n} |P(W^{(n)}<r)-e^{-\mu_{n,r}}|\le O(\frac{\ln^7 n}{n}). \lbl{th2-d1}
\end{eqnarray}

{\rm(2)} For any $a$,
\begin{equation}
 \max_{x\ge D\ln\ln n+a,C\ln n+x\in \mathbb Z}|
\exp(\frac {q}{2} p^{x}) P(W^{(n)}<C\ln n +x)-1|\le O(\frac {\ln^7 n}{n^{1-\frac {qp^a}{2}}}).\lbl{th2-d3}
\end{equation}

{\rm(3)}As $n$ tends to $\infty$,
\begin{eqnarray}
\frac {W^{(n)}-C\ln n}{\ln \ln n}\rightarrow 0 \text{ in probability }.\lbl{th2-d2}
\end{eqnarray}

{\rm(4)}  That  (\ref{d2}) holds and
\beqn
\liminf_{n\to\infty}\frac {W^{(n)}-C\ln n}{\ln\ln n}=D \,\,a.s.\lbl{th2-d4}
\eeqn

{\rm(5)} Almost surely,
\beqn
0=\liminf_{n\to\infty}\frac {W^{(2^n)}-C\ln 2^n}{\ln\ln 2^n}<\limsup_{n\to\infty}\frac {W^{(2^n)}-C\ln 2^n}{\ln\ln 2^n}=-D \lbl{th2-d5}
\eeqn
and
\beqn
\lim_{n\to\infty}\frac 1{n\ln n}(\sum_{k=1}^n W^{(2^k)}+ \frac{\ln 2}{\ln p} n^2)=0. \lbl{th2-d7}
\eeqn
\end{theorem}

Next, we shall consider the case that the success probability is not fixed.
 Suppose that $\xi^{(n)}_1,\xi^{(n)}_2,\cdots, \xi^{(n)}_n$ are i.i.d with
$P(\xi^{(n)}_i=1)=p_n=1-P(\xi^{(n)}_i=0)$. Use $U^{(n,p_n)}$ and $W^{(n,p_n)}$ to denote the
 the longest length of  arithmetic progressions or   of arithmetic progressions mod $n$  relative to $\xi^{(n)}_1,\xi^{(n)}_2,\cdots,\xi^{(n)}_n$
respectively. We have the following results.

\begin{theorem} \lbl{t3} Assume  that
\beqn
\lim_{n\to\infty} p_n=0, \lim_{n\to\infty}np_n=\infty {\text\,\, and  \,\,}
\lim_{n\to\infty} \frac {2\ln n}{-\ln p_n}=b. \lbl{th3-d0}
\eeqn

(i) If $b=\infty$, then
\beqn
\lim_{n\to\infty} P(U^{(n,p_n)}\in\{[\frac {2\ln n}{-\ln p_n}+\frac {\ln\ln n}{\ln p_n}],[\frac {2\ln n}{-\ln p_n}+\frac {\ln\ln n}{\ln p_n}]+1\})=1
\lbl{th3-d5}
\eeqn
and
\beqn
\lim_{n\to\infty} P(W^{(n,p_n)}\in \{ [\frac {2\ln n}{-\ln p_n}]-1,[\frac {2\ln n}{-\ln p_n}],[\frac {2\ln n}{-\ln p_n}]+1\}=1.\lbl{th3-d6}
\eeqn

(ii)If $b=2$, or if $2<b<\infty$ and  $b$ is not an integer, then
\beqn
\lim_{n\to\infty} P(W^{(n,p_n)}=[b])=\lim_{n\to\infty} P(U^{(n,p_n)}=[b])=1.\lbl{th3-d1}
\eeqn

(iii) If $b\ge 3$ and  $b$ is  an integer, then
\beqn
\lim_{n\to\infty} P(W^{(n,p_n)}\in \{b, b-1\})=\lim_{n\to\infty} P(U^{(n,p_n)}\in\{b,b-1\})=1.\lbl{th3-d2}
\eeqn

If in addition $u=\lim_{n\to\infty} n^2p^{b}_n\le \infty$ exists, then
\beqn
\lim_{n\to\infty} P(U^{(n,p_n)}= b-1)=e^{-\frac {u}{2(b-1)}}=1-\lim_{n\to\infty} P(U^{(n,p_n)}=b)\lbl{th3-d3}
\eeqn
and
\beqn
\lim_{n\to\infty} P(W^{(n,p_n)}=b-1)=e^{-\frac {u}{2}}=1-\lim_{n\to\infty} P(W^{(n,p_n)}=b),\lbl{th3-d4}
\eeqn
where  $[x]$  denotes the integer part of $x$.
\end{theorem}

The paper is organized as follows. In \S2, the equivalent statements of (\ref{th1-d5}) and (\ref{th2-d7}) are given.
The proofs of Theorem 1.1, Theorem 1.2 and Theorem 1.3 are given in
\S2, \S3 and \S4 respectively.

\section{Auxiliary Results}
\setcounter{equation}{0}

For clarify, we give a simple lemma that will be used.

\begin{lemma}\lbl{as-l1}  Suppose that
\beqn
b_n>0,\sum_n b_n=\infty,  \lim_{n\to\infty}\frac {\sum_{k=1}^n b_k}{n\max_{1\le k\le n} b_k}=1 \lbl{as-d1}
\eeqn
and
\beqn
a\le \liminf_{n\to\infty}  \frac {a_n}{b_n}\le \limsup_{n\to\infty} \frac {a_n}{b_n}<\infty.\lbl{as-d2}
\eeqn
 Then
\beqn
\lim_{n\to\infty} \frac {\sum_{k=1}^n a_k}{\sum_{k=1}^n b_k}=a\lbl{as-d3}
\eeqn
 if and only if for all $c>a$,
\beqn
\lim_{n\to\infty} \frac 1n \sum_{k=1}^n I_{\{{a_k}/{b_k}<c\}}=1.\lbl{as-d4}
\eeqn
\end{lemma}
\proof At first, we shall prove the  sufficiency.
 By (\ref{as-d2}), there is $d$ such that ${a_k}/{b_k}<d$ for all $k$.  For any $c>a$,
let $A_{n}=\{1\le k\le n:  {a_k}/{b_k}\ge c\}$. Then (\ref{as-d4}) implies that $\lim_{n\to\infty} {|A_{n}|}/{n}=0$.
Hence
\beq
\frac {\sum_{k=1}^n a_k}{\sum_{k=1}^n b_k}\le c+ \frac {(d-c)\sum_{k\in A_{n}}b_k}{\sum_{k=1}^n b_k}\le
c +\frac {(d-c)|A_{n}|\max_{1\le k\le n} b_k}{\sum_{k=1}^n b_k}\rightarrow c.
\eeq
The arbitrary of $c>a$ yields that  $\limsup_{n\to\infty} {\sum_{k=1}^n a_k}/{\sum_{k=1}^n b_k}\le a$.
This, together with (\ref{as-d1}) and (\ref{as-d2}), gives (\ref{as-d3}).

Next, we shall show the necessity.  By (\ref{as-d2}), for any $\varepsilon>0$, there is $K$ such that
$ {a_k}/{b_k}>a-\varepsilon$ for all $k\ge K$.  For any $c>a$, let $B_n=\{K<k\le n: {a_k}/{b_k}<c\}$.
Then for $n>K$,
\beq
\frac {\sum_{k=1}^n a_k}{\sum_{k=1}^n b_k}&\ge& c+\frac {\sum_{k=1}^{K} (a_k-cb_k)+(a-\varepsilon-c)\sum_{k\in B_n} b_k}
{\sum_{k=1}^n b_k}\\
&\ge& c+\frac {\sum_{k=1}^{K} (a_k-cb_k)+(a-\varepsilon-c)|B_{n}|\max_{1\le k\le n} b_k}{\sum_{k=1}^n b_k}.
\eeq
Letting $n$ tends to $\infty$, we get that
$$\liminf_{n\to\infty}\frac {|B_{n}|\max_{1\le k\le n} b_k}{\sum_{k=1}^n b_k}\ge \frac {c-a}{c+\varepsilon-a}.$$
The arbitrary of $\varepsilon>0$, together with (\ref{as-d1}), implies  (\ref{as-d4}) and  completes our proof.
 \qed

As an application of Lemma \ref{as-l1}. Suppose that (\ref{as-d2}) holds. In addition, suppose that
$b_nc_n>0$, $\sum_{n} b_nc_n=\infty$ and $ \lim_{n\to\infty}\frac {\sum_{k=1}^n b_kc_k}{n\max_{1\le k\le n} b_kc_k}=1$.
Then $\lim_{n\to\infty}  {\sum_{k=1}^n  a_kc_k}/{\sum_{k=1}^n b_kc_k}=a$ if and only if
(\ref{as-d4}) holds for all $c>a$.
Particularly, by letting $c_n=1/b_n$,   we see that   $\lim_{n\to\infty} \frac 1n  \sum_{k=1}^n  \frac {a_k}{b_k}=a$
if and only if (\ref{as-d4}) holds for all $c>a$.
Therefore
 if (\ref{as-d1}) and (\ref{as-d2}) hold,  then
 $\lim_{n\to\infty} {\sum_{k=1}^n  a_k}/{\sum_{k=1}^n b_k}=a$ if and only if
 $\lim_{n\to\infty} \frac 1n  \sum_{k=1}^n  \frac {a_k}{b_k}=a$, if and only if (\ref{as-d4}) holds for all $c>a$.

Note that $n\ln n-n\le \sum_{k=1}^n \ln k \le n\ln n$. By Lemma \ref{as-l1},
we have the following propositions.

\begin{prop}\lbl{as-p1}If (\ref{th1-d4}) holds, then  (\ref{th1-d5}) holds if and only if
\beqn
\lim_{n\to\infty}\frac 1n {\sum_{k=2}^n \frac{U^{(2^k)}-C\ln 2^k}{\ln\ln 2^k}}=D \,\,a.s, \lbl{as-d5}
\eeqn
and also  if and only if for any $1>\varepsilon>0$,
\beqn
\lim_{n\to\infty}\frac 1n \sum_{k=1}^n I_{\{U^{(2^k)}<C\ln 2^k+D(1-\varepsilon)\ln\ln 2^k\}}=1 \,\, a.s. \lbl{as-d6}
\eeqn
\end{prop}

\begin{prop}\lbl{as-p3} If (\ref{th2-d5}) holds, then  (\ref{th2-d7}) holds if and only if
\beqn
\lim_{n\to\infty}\frac 1n {\sum_{k=2}^n \frac{W^{(2^k)}-C\ln 2^k}{\ln\ln 2^k}}=0 \,\,a.s, \lbl{as-d9}
\eeqn
and also if and only if for any $\varepsilon>0$,
\beqn
\lim_{n\to\infty}\frac 1n \sum_{k=1}^n I_{\{W^{(2^k)}<C\ln 2^k-\varepsilon D\ln\ln 2^k\}}=1 \,\, a.s. \lbl{as-d10}
\eeqn
\end{prop}

\section{The asymptotic distribution of $U^{(n)}$}

\setcounter{equation}{0}

Suppose that  $2\le r\le n$. Let
\begin{eqnarray*}
B_{n}&=&B^{(r)}_{n}=\{(a,s):1\le a,s\le n,  a+(r-1)s\le n\}\\
&=&\{(a,s):1\le s\le [\frac {n-1}{r-1}],1\le a\le n-(r-1)s\}.
\end{eqnarray*}
For any $(a,s)\in B_n$, let
 $$A_{a,s}=A^{(r)}_{a,s}=\{\xi_a=1,\xi_{a+s}=1,\cdots, \xi_{a+(r-1)s}=1\}\cap\{a-s\le 0 \text{ or }\xi_{a-s}=0\}.$$
 Then
  $$P(U^{(n)}\ge r)=P(\bigcup\limits_{(a,s)\in B_n} A_{a,s}).$$
  Let $I=I_{n,r}=\sum\limits_{(a,s)\in B_n} P(A_{a,s})$. Then we have
\begin{eqnarray}
P(U^{(n)}\ge r)\le I. \lbl{s2-d1}
\end{eqnarray}
 Set
 $$B_{a,s}=B^{(r)}_{a,s}=\left\{
             \begin{array}{ll}
               \{a,a+s,\cdots, a+(r-1)s\}, & \hbox{if $a\le s$;} \\
               \{a-s,a,a+s,\cdots, a+(r-1)s\}, & \hbox{otherwise.}
             \end{array}
           \right.$$
Let $G$ be the graph with vertex set $B_n$ and edges defined by
$(a,s) \thicksim (b,t)$ if and only if $B_{a,s}\cap B_{b,t}\not=\emptyset$.
Then $G$ is a dependency graph of $\{I_{A_{a,s}}: (a,s)\in B_n\}$, where $I_{A_{a,s}}$ is the indicator
function of $A_{a,s}$.
 The notion of dependency graphs can be found in \cite{AS00},\cite{BC05} or in \S2.1 of \cite{BYZ07}.
Set
$$e^{(n,r)}=\sum_{(a,s)\in B_n}\sum_{(b,t)\thicksim (a,s)} P(A_{a,s})P(A_{b,t})+ \sum_{(a,s)\in B_n}\sum_{(b,t)\thicksim (a,s),(b,t)\not= (a,s)} P(A_{a,s}A_{b,t}).$$
Note that $P(U^{(n)}<r)=P(\sum_{(a,s)\in B_n} I_{A_{a,s}}=0)$.
Applying  the Chen-Stein method, (see \cite{AGG89}, \cite{BC05} or Theorem 3 of  \cite{BYZ07}),  we get that
\begin{eqnarray}
|P(U^{(n)}<r)-e^{-I}|\le e^{(n,r)}. \lbl{s2-d2}
\end{eqnarray}

 \begin{lemma} \lbl{s2-l1}
It holds that
\begin{eqnarray}
p^r\frac {(n-r)^2}{2(r-1)}-p^{r+1}\frac {n^2}{2r}\le I\le p^r\frac {n^2}{2(r-1)}-p^{r+1}\frac {(n-r)^2}{2r}.\lbl{s2-l1-d1}
\end{eqnarray}
 \end{lemma}
  \proof
Clearly,   $|B_n|=\frac {[\frac {n-1}{r-1}]}{2} (2n-r+1-(r-1)[\frac {n-1}{r-1}])$.
It implies that
\beqn
\frac {(n-r)^2}{2(r-1)}\le|B_n|\le \frac {n^2}{2(r-1)}. \lbl{s2-l1-d2}
\eeqn
Note that
$B_n\cap \{(a,s):a>s\}=\{(a,s):1\le s\le [ {n-1\over r}], s<a\le n-(r-1)s\}$. We have
\beqn
\frac {(n-r)^2}{2r}\le |B_n\cap \{(a,s):a>s\}|\le \frac {n^2}{2r}. \lbl{s2-l1-d3}
\eeqn
Clearly,
$I=p^r|B_n|-p^{r+1}|B_n\cap\{(a,s):a>s\}|$.
This, together with (\ref{s2-l1-d2}) and (\ref{s2-l1-d3}), gives (\ref{s2-l1-d1}).
\qed

\begin{lemma} \lbl{s2-l2}It holds that
 $e^{(n,r)}\le  9(n^3 p^{2r-1}+n^2 r^3 p^{\frac 53 r-1}+n^2 p^{\frac 32 r-1})$.
\end{lemma}
\proof
Let
\begin{eqnarray*}
c_1&=&|\{(a,s,b,t)\in B_n\times B_n:  (a,s)\thicksim(b,t)\}|,\\
c_2&=&|\{(a,s,b,t)\in B_n\times B_n: |B_{a,s}\cap B_{b,t}|\ge 2\}|
\end{eqnarray*}
and
\beq
c_3=|\{(a,s,b,t)\in B_n\times B_n: (a,s)\thicksim(b,t),t=2s\text{ or } s=2t\}|.
\end{eqnarray*}
Then
\begin{eqnarray*}
c_1\le &&|\{(a,s,b,t): (a,s)\in B_n, 1\le t\le [\frac {n-1}{r-1}], b=a+is-jt,\\
&&-1\le i,j\le r-1\}|
\le 9n^3/2
\end{eqnarray*}
and
\begin{eqnarray*}
c_3&\le& 2|\{(a,s,b,t): (a,s)\in B_n, t=2s , b=a+is-jt,-1\le i,j\le r-1\}|\\
&=&2|\{(a,s,b):(a,s)\in B_n, b=a+ks,-2r+1\le k\le r+1\}|\le 7n^2.
\end{eqnarray*}

Suppose that $|B_{a,s}\cap B_{b,t}|\ge 2$ and $x_0$ is the minimal number of the set
$B_{a,s}\cap B_{b,t}$.  Then $x_0=a+is=b+jt$ for some $-1\le i,j\le r-1$.
If $x\in B_{a,s}\cap B_{b,t}$ and $x>x_0$, then $x=a+i's=b+j't$ for some $-1\le i',j'\le r-1$.
It follows that $x-x_0=(i'-i)s=(j'-j)t$. Thus $t=\frac {k_1}{k_2}s$ for some $1\le k_1,k_2 \le r$.
In addition,  there is  positive integer $k$ such that
$i'-i=kt_0, j'-j=ks_0$ and
$x-x_0=kst_0$, where $s_0={s}/{\gcd(s,t)}$ and
$t_0={t}/{\gcd(s,t)}$.
Since
$i'-i\le r$, $k\le {r}/{t_0}$. Similarly, $k\le  {r}/{s_0}$. Therefore
$$|B_{a,s}\cap B_{b,t}|\le  {r}/\max(s_0, t_0)+1.$$
Consequently,  $|B_{a,s}\cap B_{b,t}|\le {r}/{3}+1$ whenever $\max(s_0, t_0)\ge 3$.
When   $\max(s_0, t_0)=2$,   $|B_{a,s}\cap B_{b,t}|\le {r}/{2}+1$. Actually, in this case, $s=2t$ or $t=2s$.
When $\max(s_0, t_0)=1$,  $s=t$.
We shall show that $A_{a,s}\cap A_{b,s}=\emptyset$ whenever $a\not=b$ and $B_{a,s}\cap B_{b,s}\not =\emptyset$.
Assume that $b>a$ without loss of generality. Since $B_{a,s}\cap B_{b,s}\not =\emptyset$,  $a+is=b+js$ for some $-1\le i,j\le r-1$ and
 hence $b=a+ks$ for some  $1\le k \le r$.
Thus $A_{a,s}\subseteq\{\xi_{a+(k-1)s}=1\}$ and $A_{b,s}\subseteq\{\xi_{b-s}=0\}=\{\xi_{a+(k-1)s}=0\}$.
It implies that $A_{a,s}\cap A_{b,s}=\emptyset$ as desired.
In view of the discussion above, we have
\begin{eqnarray*}
c_2\le&& |\{(a,s,b,t): (a,s)\in B_n, t= {sk_1}/{k_2}, b=a+is-jt,\\
&&-1\le i,j\le r-1,1\le k_1,k_2\le r\}|
\le  9n^2r^3/4
\end{eqnarray*}
and
\begin{eqnarray*}
e^{(n,r)}
\le 2c_1p^{2r-1}+c_2p^{\frac {5r}{3}-1}+c_3p^{\frac {3r}{2}-1}
\le &9(n^3 p^{2r-1}+n^2 r^3 p^{\frac 53 r-1}+n^2 p^{\frac 32 r-1})
\end{eqnarray*}
as desired.
\qed

Similar as the proof of Lemma \ref{s2-l2}, we may show that for any $m$, $n$ and $2\le r_m\le r_n$, with
$H=\{(a,s,b,t):(a,s)\in B^{(r_m)}_m, (b,t)\in B^{(r_n)}_n\}$,
\beqn
|\{(a,s,b,t)\in H:|B^{(r_m)}_{a,s}\cap B^{(r_n)}_{b,t}|\ge 1\}|\le9m^2n/2 \lbl{s2-l3-d1}
\eeqn
and
\beqn
|\{(a,s,b,t)\in H:|B^{(r_m)}_{a,s}\cap B^{(r_n)}_{b,t}|\ge 2\}|\le  9m^2r_mr^2_n/4.\lbl{s2-l3-d2}
\eeqn
Also we can show that if $|B^{(r_m)}_{a,s}\cap B^{(r_n)}_{b,t}|>\frac {r_n}{2}+1$ and
$A^{(r_m)}_{a,s}\cap A^{(r_n)}_{b,t}\not =\emptyset$, then $(b,t)=(a,s)$.
Therefore
\beqn
|\{(a,s,b,t)\in H:|B^{(r_k)}_{a,s}\cap B^{(r_m)}_{b,t}|>\frac {r_n}{2}+1,A^{(r_m)}_{a,s}\cap A^{(r_n)}_{b,t}\not =\emptyset\}|\le  k^2/r_k.\lbl{s2-l3-d3}
\eeqn

 \textbf{Proof of Theorem 1.1}
(1) Lemma \ref{s2-l1} and   (\ref{s1-d1}) imply that
\beqn
|e^{-I_{n,r}}-e^{-\lambda_{n,r}}|\le |I_{n,r}- \lambda_{n,r}|\le 2np^r.\lbl{s2-pd1}
\eeqn
 Let $r_n=[\frac {-2\ln n}{\ln p}+2 \frac {\ln\ln n}{\ln p}+\frac {\ln\ln\ln n}{2\ln p}]$
and
$R_n=[\frac {-3\ln n}{\ln p}]$.
By (\ref{s2-d2}), (\ref{s2-pd1}) and Lemma \ref{s2-l2},
\begin{eqnarray}
&&\max_{r_n\le r\le R_n}|P(U^{(n)}<r)-e^{-\lambda_{n,r}}|\nn\\
&&\le \max_{r_n\le r\le R_n}\big(|P(U^{(n)}<r)-e^{-I_{n,r}}|+|e^{-I_{n,r}}-e^{-\lambda_{n,r}}|\big)\nn\\
&&\le  \max_{r_n\le r\le R_n} (e^{(n,r)}+ 2np^{r})\nn\\
&&\le 9(n^3 p^{2r_n-1}+n^2 R^3_n p^{\frac 53 r_n-1}+n^2 p^{\frac 32 r_n-1})+2np^{r_n}\nn\\
&&=O(\frac {\ln^{4} n\ln\ln n}{n}).\lbl{th1-pd100}
\end{eqnarray}
On the other hand, it's easy to check that
$e^{-\lambda_{n,r_n}}=e^{-O(\ln n\sqrt{\ln\ln n})}=o(n^{-1})$
and
$1-e^{-\lambda_{n,R_n}}=1-e^{-O(\frac 1{n\ln n})}=o(n^{-1})$.
Note that
$P(U^{(n)}<r)$ and $e^{-\lambda_{n,r}}$ are both increasing functions of $r$.
Hence when $r<r_n$,
\beq
&&|P(U^{(n)}<r)-e^{-\lambda_{n,r}}|\le P(U^{(n)}<r)+e^{-\lambda_{n,r}}\le P(U^{(n)}<r_n)+e^{-\lambda_{n,r_n}}\\
&\le& |P(U^{(n)}<r_n)-e^{-\lambda_{n,r_n}}|+2e^{-\lambda_{n,r_n}}\le O(\frac {\ln^{4} n\ln\ln n}{n}).
\eeq
Similarly, when $r>R_n$,
\beq
&&|P(U^{(n)}<r)-e^{-\lambda_{n,r}}|\le 1-P(U^{(n)}<R_n)+1-e^{-\lambda_{n,R_n}}\\
&\le& |P(U^{(n)}<R_n)-e^{-\lambda_{n,R_n}}|+2(1-e^{-\lambda_{n,R_n}})\le O(\frac {\ln^{4} n\ln\ln n}{n}).
\eeq
This completes the proof of (\ref{th1-d1}).\qed

(2) Let $r=h_n+x$. For convenience, set
\begin{eqnarray*}
\epsilon_{n,x}=|\frac {-q\ln p}{4} p^{x}-\lambda_{n,r}|= p^{x}\big|\frac {-q\ln p}{4}-\frac {q\ln n}{2(r-1)}-\frac{p \ln n}{2r(r-1)}\big|.
\end{eqnarray*}
Then
\begin{eqnarray}
&&|\exp(\frac {-q\ln p}{4} p^{x}) P(U^{(n)}<h_n +x)-1|\nn\\
&\le & \exp(\frac {-q\ln p}{4} p^{x})|P(U^{(n)}<r)-e^{-\lambda_{n,r}}|+|\exp(\epsilon_{n,x})-1|\lbl{s2-pd2}
\end{eqnarray}
Since $a<1$,
 $$\exp(\frac {-q\ln p}{4} p^{aD\ln\ln n})=\exp(\frac {-q\ln p}{4}\ln^a n)=o(n^{\frac 13}).$$
Thus
\begin{eqnarray}
&&\max_{x\ge aD\ln\ln n,h_n+x\in \mathbb Z}\exp(\frac {-q\ln p}{4} p^{x})|P(U^{(n)}<r)-e^{-\lambda_{n,r}}|\nn\\
&\le & o(n^{\frac 13})O(n^{-1}\ln^{4} n\ln\ln n))=o(n^{-\frac 12}).\lbl{s2-pd3}
\end{eqnarray}
It's easy to verify that
\begin{eqnarray}
\max_{  x > -D\ln\ln n+1}\epsilon_{n,x}\le
p^{-D\ln\ln n}O(1)=O(\frac 1{\ln n})\lbl{s2-pd4}
\end{eqnarray}
and
\begin{eqnarray}
\max_{ aD\ln\ln n\le x \le -D\ln\ln n+1}\epsilon_{n,x}\le
p^{ aD\ln\ln n}O(\frac {\ln\ln n}{\ln n})=O(\frac{ \ln\ln n}{\ln^{1-a} n}).\lbl{s2-pd5}
\end{eqnarray}
Now (\ref{th1-d3}) follows  by (\ref{s2-pd2})--(\ref{s2-pd5}).

(3) For any  $\varepsilon>0$, (\ref{s2-d1}) and (\ref{s2-l1-d1}) imply that
\begin{eqnarray}
P(U^{(n)}\ge C\ln n +(1-\varepsilon)D\ln\ln n)\le O(\ln^{-\varepsilon} n).\lbl{th1-pd7}
\end{eqnarray}
On the other hand,  by (\ref{s1-d1}) and (\ref{th1-d1}),
\begin{eqnarray}
P(U^{(n)}< C\ln n +(1+\varepsilon)D\ln\ln n)\le e^{-O(\ln^{\varepsilon} n)}+O(\frac {\ln^{4} n\ln\ln n}{n}).\lbl{th1-pd8}
\end{eqnarray}
Hence (\ref{th1-d2}) holds.

(4) By (\ref{th1-pd8}),
$\sum_{k=1}^\infty P(U^{(2^k)}< C\ln 2^k +(1+\varepsilon)D\ln\ln 2^k)<\infty $
for any $\varepsilon>0$.
 One then deduces from the Borel-Cantelli Lemma that  $$P(U^{(2^k)}<C\ln 2^k +(1+\varepsilon)D \ln\ln 2^k \,\,\,\,i.o.)=0.$$ It follows that for almost surely $\omega$, there is $K(\omega)$ such that for $k\ge K(\omega)$,
\beqn
U^{(2^k)}(\omega)\ge C\ln 2^k +(1+\varepsilon)D\ln\ln 2^k.\lbl{th1-pd9}
\eeqn
If $n>2^{K(\omega)}$, then $2^k\le n<2^{k+1}$ for some $k\ge K(\omega)$. Hence
$U^{(2^k)}(\omega)\le U^{(n)}(\omega)\le U^{(2^{k+1})}(\omega)$. This, together with (\ref{th1-pd9}), gives that
$$U^{(n)}(\omega)\ge  C\ln n-C\ln 2+(1+\varepsilon)D\ln\ln n.$$
Now the arbitrary of $\varepsilon>0$ yields that $\liminf_{n\to\infty} \frac {U^{(n)}-C\ln n}{\ln\ln n}\ge D$.
Therefore  $\liminf_{n\to\infty} \frac {U^{(n)}-C\ln n}{\ln\ln n}= D$   by considering (\ref{th1-d2}).

Let $T_k$ be the longest length of  arithmetic progressions relative to
$\xi_{2^{k-1}+1},\\ \xi_{2^{k-1}+2},\cdots,\xi_{2^{k}}$.
Then $T_1,T_2,\cdots$ are independent. In addition,  $T_k$ has the same distribution with  $U^{(2^{k-1})}$. By
(\ref{th1-d1}), $$P(T_k\ge C\ln 2^k)=P(U^{(2^{k-1})}\ge C\ln 2^k)\le O(1/k).$$
Hence $\sum_{k=1}^\infty P(T_k\ge C\ln 2^k)=\infty$. By Borel-Cantelli Lemma,
$P(T_k\ge C\ln 2^k\,\,\, i.o.)=1$.
Consequently,  $\limsup_{k\to\infty} \frac {U^{(2^k)}-C\ln 2^k}{\ln\ln 2^k}\ge 0$ by noting that $U^{(2^k)}\ge T_k$.
On the other hand,
 (\ref{th1-pd7}) yields that
$$\sum_{k=1}^\infty P(U^{(2^k)}\ge C\ln 2^k +(1-\varepsilon)D \ln\ln 2^k)<\infty $$
whenever  $\varepsilon>1$. Hence $\limsup_{k\to\infty} \frac {U^{(2^k)}-C\ln 2^k}{\ln\ln 2^k}\le 0$.
 Therefore
 $\limsup_{k\to\infty} \\\frac {U^{(2^k)}-C\ln 2^k}{\ln\ln 2^k}= 0$. Furthermore, we can deduce that
$\limsup_{n\to\infty} \frac {U^{(n)}-C\ln n}{\ln\ln n}= 0$ by the fact that
$U^{(n)}$ is increasing. This completes the proof of (\ref{th1-d4}).

Now we come to  prove (\ref{th1-d5}).
By Proposition \ref{as-p1}, we need only to show (\ref{as-d6}).
Let $r_k= [C\ln 2^k+D(1-\varepsilon)\ln\ln 2^k]$, $V_n=\{(a,s,k):1\le k\le n, (a,s)\in B_{2^k}^{(r_k)}\}$ and
$\Lambda(n)=\sum_{(a,s,k)\in V_n}I_{ A_{a,s}^{(r_k)}}$.
Then it suffices to show that $\lim_{n\to\infty}\Lambda(n)/n=0$ a.s.
Clearly, $E\Lambda(n)=\sum_{k=1}^n I_{2^k,r_k}=O(n^{1-\varepsilon})$ and
$D\Lambda (n)$ is less than the sum of $p^{r_k+r_m-|B^{(r_k)}_{a,s}\cap B^{(r_m)}_{b,t}|}$ with
$(a,s,k),(b,t,m)\in V_n$, $B^{(r_k)}_{a,s}\cap B^{(r_m)}_{b,t}\not=\emptyset$ and $A^{(r_k)}_{a,s}\cap A^{(r_m)}_{b,t}\not=\emptyset$.
By (\ref{s2-l3-d1})--(\ref{s2-l3-d3}),
\beq
D\Lambda(n)\le &&\sum_{1\le i\le j\le n} O(2^{-j}ij+2^{-j}i^{2}j^3+2^{2i-2j}j^{1-\varepsilon}i^{-1})\\
\le&& O(\sum_{j=1}^\infty 2^{-j} j^{6})+O(\sum_{j=1}^n j^{-\varepsilon}\sum_{k=0}^{j-1} 2^{-2k}j(j-k)^{-1})\\
\le &&O(1)+\sum_{k=0}^{\infty} 2^{-2k}(k+1)O(\sum_{j=1}^n j^{-\varepsilon})=O(n^{1-\varepsilon}).
\eeq
Then by Tchebychev's inequality, for any $\delta>0$,
$$\sum_{n=1}^\infty P(|\Lambda_n/n-E(\Lambda_n/n)|>\delta)\le \sum_{n=1}^\infty O(\frac 1{n^{1+\varepsilon}\delta^2})<\infty.$$
 The Borel-Cantelli Lemma  yields that  $\Lambda_n/n\to 0\,\,a.s.$. Hence
  (\ref{as-d6}) holds as desired.

Finally, we shall prove  (\ref{th1-d6}).
Let $c_n=[\ln n/\ln 2]$. Then $2^{c_n}\le n<2^{c_n+1}$. For any integers $1\le a\le b$,
\beq
\ln\frac{b+1}a=\int_{a}^{b+1} \frac 1x\,dx\le\sum_{i=a}^b\frac 1i\le \int_{a-1}^b \frac 1x\,dx=\ln \frac b{a-1}.
\eeq
Thus
$$\sum_{k=1}^n \frac {U^{(k)}}{k}\ge\sum_{i=0}^{c_n-1} U^{(2^i)}\sum_{j=2^i}^{2^{i+1}-1}\frac 1j \ge \ln 2\sum_{i=0}^{c_n-1} U^{(2^i)}$$
and
$$\sum_{k=2}^n \frac {U^{(k)}}{k}\le\sum_{i=0}^{c_n} U^{(2^{i+1})}\sum_{j=2^i+1}^{2^{i+1}}\frac 1j \le \ln 2\sum_{i=1}^{c_n+1} U^{(2^i)}.$$
Hence by (\ref{th1-d5}),
$\lim_{n\to\infty}\frac 1{c_n\ln c_n}(\frac 1{\ln 2}\sum_{k=1}^n \frac {U^{(k)}}{k}+\frac{\ln 2}{\ln p} c^2_n)=D\,\,a.s.$
It follows (\ref{th1-d6}) immediately and completes the proof of Theorem 1.1.
\qed

\section{The asymptotic distribution of $W^{(n)}$}
\setcounter{equation}{0}
Suppose that $2\le r\le n$. For any $1\le a,s \le n$, let
\beqn
 \tilde A_{a,s}=\tilde A^{(n,r)}_{a,s}=\{\xi_a=0,\xi_{a+s \mod n}=1,\cdots, \xi_{a+rs \mod n}=1\}.\lbl{s3-dd1}
\eeqn
Let
$$\tilde B_n=\tilde B^{(r)}_n=\{(a,s): 1\le a\le n,1\le s\le [ n/2],\gcd(n,s)<n/r\}$$
 and
\begin{eqnarray}
A_1=\bigcup\limits_{(a,s)\in \tilde B_n} \tilde A_{a,s}. \lbl{s3-d1}
\end{eqnarray}
Set $C_{a,s}=\{a,a+s \mod n,a+2s \mod n,a+3s \mod n,\cdots\}$.  Then
$C_{a,s}=\{a,a+s \mod n,\cdots, a+ (\frac {n}{\gcd(s,n)}-1)s\mod n \}$  and $|C_{a,s}|=n/\gcd(s,n)$. Set
\beqn
A_2=\bigcup\limits_{s|n,s\le n/r,1\le a\le s} \{\xi_i=1,\forall i\in C_{a,s}\},\lbl{s3-d2}
\eeqn
where $s|n$ means that $s$ is a divisor of $n$.

\begin{lemma}\lbl{s3-l1} It holds that
\beqn
\{W^{(n)}\ge r\}=A_1\cup A_2.\lbl{s3-l1-d1}
\eeqn
\end{lemma}
\proof
Put $$W^{(n)}_s=\max\limits_{1\le i \le n} W^{(n)}_{i,s},$$
which is the maximum length of  arithmetic progressions mod $n$ in $\Sigma_n$ with difference $s$.
For any $m \ge 0$,  $\{\xi_a=1,\xi_{a+s \mod n}=1,\cdots,\xi_{a+ms \mod n}=1\} $  if and only if
$\{\xi_b=1,\xi_{b+(n-s) \mod n}=1,\cdots,\xi_{b+m(n-s) \mod s}=1\} $, where $b=a+ms \mod n$. In addition,
$\gcd(s,n)=\gcd(n-s,n)$.
Hence $W_{s}^{(n)}=W_{n-s}^{(n)}$ for all $1\le s \le n$. Consequently,
$$W^{(n)}=\max_{1\le s \le n} W^{(n)}_s=\max_{1\le s \le [ n/2]} W^{(n)}_s.$$

For any $1\le a,b\le n$, $C_{a,s}\cap C_{b,s}=\emptyset$ or  $C_{a,s}=C_{b,s}$.
In addition $C_{a,s}=C_{b,s}$ if and only if $b=a+\gcd(s,n)\cdot k$ for some integer $k$.
Thus $\{1,2,\cdots,n\}$ is the
disjoint union of $C_{a,s}$ with $1\le a\le \gcd(s,n)$. It follows that
\beq
W^{(n)}_s=\max_{1\le a \le \gcd(s,n)}\tilde W^{(n)}_{a,s}.
\eeq
where $\tilde W^{(n)}_{a,s}=\max\limits_{i\in C_{a,s}}  W^{(n)}_{i,s}$.
Note that
$\{\tilde W^{(n)}_{a,s}\ge r\}= \{\xi_i=1,\forall i\in C_{a,s}\}$ when $n/\gcd (s,n)= r$,
and
\beq
\{\tilde W^{(n)}_{a,s}\ge r\}=(\bigcup\limits_{i\in C_{a,s}} \tilde A_{i,s})\bigcup \{\xi_i=1,\forall i\in C_{a,s}\}.
\eeq
provided $n/\gcd (s,n)> r$.
Therefore
\begin{eqnarray*}
\{W^{(n)}\ge r\}= \big(\bigcup_{(i,s)\in\tilde B_n}\tilde A_{i,s}\big)\bigcup
\big(\bigcup_{1\le a\le \gcd(s,n),1\le s \le [ \frac n2],  \frac n{\gcd(s,n)}\ge r}\{\xi_i=1,\forall i\in C_{a,s}\}\big).
\end{eqnarray*}
This, together with  the fact that
$C_{a,s}=C_{a,\gcd(s,n)}$,  yields (\ref{s3-l1-d1}). \qed

Let $\tilde I=\tilde I_{n,r}=\sum_{(a,s)\in \tilde B_n} P(\tilde A_{a,s})$.
By Lemma (\ref{s3-l1}), we have
\beqn
 P(A_1)\le P(W^{(n)}\ge r)\le P(A_1)+P(A_2)\le \tilde I+P(A_2).\lbl{s3-d3}
\eeqn
Set
\beqn
\tilde B_{a,s}=\tilde B^{(n,r)}_{a,s}=\{a,a+s \mod n,\cdots, a+rs \mod n\}.\lbl{s3-dd2}
\eeqn
Define  $\tilde G$ to be the graph with vertex set $\tilde B_n$ and edges defined by
$(a,s) \thicksim (b,t)$ if and only if $\tilde B_{a,s}\cap \tilde B_{b,t}\not=\emptyset$.
Then $\tilde G$ is a dependency graph of $\{I_{\tilde A_{a,s}}: (a,s)\in \tilde B_n\}$.
Put
$$\tilde e^{(n,r)}=\sum_{(a,s)\in \tilde B_n}\sum_{(b,t)\thicksim (a,s)} P(\tilde A_{a,s})P(\tilde A_{b,t})+\sum_{(a,s)\in \tilde B_n}\sum_{(b,t)\thicksim (a,s),(b,t)\not= (a,s)} P(\tilde A_{a,s} \tilde A_{b,t}).$$
Then
\beqn
|P(A^c_1)-e^{-\tilde I}|\le \tilde e^{(n,r)}.\lbl{s3-d4}
\eeqn
The estimations of $P(A_2)$, $\tilde I$ and $\tilde e^{(n,r)}$ are given in the following two lemmas.

\begin{lemma}\lbl{s3-l2} We have
\beqn
P(A_2)\le  \frac {np^r}{qr}\lbl{s3-l2-d1}
\eeqn
and
\beqn
(1 -\frac {(r+1)^2}{2n})\frac {qn^2p^r}{2}\le \tilde I \le \frac {qn^2 p^r}{2}.\lbl{s3-l2-d2}
\eeqn
\end{lemma}
\proof Obviously,
\begin{eqnarray*}
P(A_2)&\le& \sum_{s|n,s\le  n/r}sp^{\frac ns}\le \sum_{i=r}^n  \frac nip^i\le \frac nr \sum_{i=r}^n p^i\le \frac {np^r}{qr}.
\end{eqnarray*}
Clearly,
$\{1\le s\le [\frac n2]: \gcd (n,s)\ge \frac nr\}\subseteq \{s=\frac ni j:2\le i\le r,1\le j\le [\frac i2]\}$.
It implies that
\begin{eqnarray*}
\frac {n^2}{2}\ge |B_n|\ge n([\frac n2]-\sum_{i=2}^r [\frac i2])\ge \frac {n^2}2(1 -\frac {(r+1)^2}{2n}).
\end{eqnarray*}
Now the fact that
$\tilde I=|\tilde B_n|qp^r$ yields (\ref{s3-l2-d2}) immediately.
\qed

\begin{lemma}\lbl{s3-l3}It holds that $\tilde e^{(n,r)}\le 4(n^3r^2p^{2r-1}+n^2r^5p^{\frac 32  r-1}+nr^6 p^r)$.
\end{lemma}
\proof
Let $H=\{(a,s,b,t)\in \tilde B_n\times \tilde B_n: (a,s)\sim (b,t)\}$, $\tilde c_1=|H|$,
\begin{eqnarray*}
\tilde c_2&=&|\{(a,s,b,t)\in H: |\tilde B_{a,s}\cap \tilde B_{b,t}|\ge 2\}|
\end{eqnarray*}
and
\begin{eqnarray*}
\tilde c_3=|\{(a,s,b,t)\in H: (a,s)\not=(b,t),|\tilde B_{a,s}\cap \tilde B_{b,t}|>\frac r 2+1,\tilde A_{a,s}\cap\tilde A_{b,t}\not=\emptyset\}|.
\end{eqnarray*}
Then
\begin{eqnarray}
\tilde c_1\le&& |\{(a,s,b,t): (a,s)\in \tilde  B_n, 1\le t\le [\frac n2], b=a+is-jt \mod n,\nn\\
&&0\le i,j\le r\}|
\le n^3r^2.\lbl{s3-l3-d1}
\end{eqnarray}
Suppose that $|\tilde B_{a,s}\cap \tilde B_{b,t}|\ge 2$.
Then there is $0\le j_1<j_2\le r$ and $x,y\in \tilde B_{a,s}$ such that $b+j_1t\mod n=x$ and $b+j_2t\mod n=y$.
Hence $(j_2-j_1)t-kn=y-x$ for some $0\le k \le j_2-j_1$ and $b=x-j_1t\mod n$.  Therefore
 \begin{eqnarray}
\tilde c_2\le&& |\{(a,s,b,t): (a,s)\in \tilde  B_n, t=(kn+y-x)/i, b=x-jt\mod n,\nn\\
&&1\le i\le r, 0\le j,k\le r,x,y\in \tilde B_{a,s}\}|
\le 3n^2r^5.\lbl{s3-l3-d1}
\end{eqnarray}

We shall show that $\tilde A_{a,s}\cap\tilde A_{b,s}=\emptyset$ when   $b\not=a$  and $(a,s)\thicksim (b,s)$.
Let $$j^*=\min\{0 \le j \le r: b+j s \mod n \in \tilde B_{a,s}\}.$$
Then there is $0\le i \le r$ such that
$b+j^*  s \mod n =a+ is \mod n$.
It follows that  $b+(j^*-1)  s \mod n =a+ (i-1)s \mod n$. Hence $j^*-1<0$ or $i-1<0$.
If $j^*-1<0$, then $j^*=0$ and hence $b=a+is \mod n$. Since $b\not =a$, $i\not=0$.
We have
$\tilde A_{a,s}\cap\tilde A_{b,s}\subseteq\{\xi_{a+is \mod n}=1,\xi_{b}=0\}=\{\xi_{b}=1,\xi_{b}=0\}=\emptyset$.
Similarly  $\tilde A_{a,s}\cap\tilde A_{b,s}=\emptyset$ when $i-1<0$.

Suppose that $s\not=t$, $|\tilde B_{a,s}\cap \tilde B_{b,t}|> r/2+1$ and
$|\tilde B_{a,s}\cap \tilde B_{b,t}|=\{a+i_0s,a+i_1s,\cdots,a+i_ks\}$ with
$0\le i_0<i_1<\cdots<i_k\le r$. Then there is $ l$ such that
$i_{l+1}-i_l=1$. Since  $a+i_ls \mod n=b+j_1 \mod n$ and $a+i_{l+1}s\mod n =b+j_2 \mod n$ for some
$0\le j_1,j_2\le r$ and $j_1\not=j_2$,  $s=it \mod n$ with $i=j_2-j_1$.
If $i=1$, then $s=t$. If $i=-1$, then $s=n-t$ and hence $s=t=n/2$ by the fact that
$1\le s,t\le n/2$. The contradiction shows that $1<|i|\le r$.
Similarly,
$t=js \mod n$ for some $1< |j|\le r$.  It follows that
$s=ijs \mod n$, that is $(ij-1)s=\jmath n$ for some $|\jmath|\le r^2/2$.
Consequently,
\begin{eqnarray}
\tilde c_3 &\le& |\{(a,s,b,t): s=\jmath n/(ij-1), t=js \mod n,  b=a+ls-mt\mod n,\nn\\
&&1\le a\le n, 1<|i|,|j|\le r,  |\jmath|\le r^2/2, 0\le l,m\le r\}|\le  4nr^6. \lbl{s3-l3-d6}
\end{eqnarray}
Thus our result holds by noting that
$\tilde e^{(n,r)}\le 2\tilde c_1 p^{2r-1}+\tilde c_2 p^{\frac{3r}{2}-1}+\tilde c_3 p^{r}$.
\qed

Similarly,  we may show that for any fixed $m, n$,
\beqn
|\{(a,s,b,t)\in \tilde B_m^{(r_m)}\times \tilde B_n^{(r_n)}:|\tilde B^{(m,r_m)}_{a,s}\cap \tilde B^{(n,r_n)}_{b,t}|\ge 1\}|\le m^2nr_mr_n \lbl{s3-l4-d1}
\eeqn
and
\beqn
|\{(a,s,b,t)\in  \tilde B_m^{(r_m)}\times \tilde B_n^{(r_n)}:|\tilde B^{(m,r_m)}_{a,s}\cap \tilde B^{(n,r_n)}_{b,t}|\ge 2\}|\le  3m^2r^2_mr^3_n. \lbl{s3-l4-d2}
\eeqn

\begin{lemma}\lbl{s3-l4}
Suppose that  $(a,s)\in \tilde B_{m}^{(r_m)}$ and $(b,t)\in \tilde B_{n}^{(r_n)}$.
If $n\ge 2m$, $r_n\ge 36$ and $t>3n/r_n$,
then
$|\tilde B^{(m,r_m)}_{a,s}\cap \tilde B^{(n,r_n)}_{b,t}|<\frac 34 r_n$.
\end{lemma}
\proof If $1\le x\le m$, $A=\{x,x+t,\cdots,x+kt\}\subseteq\{1,2,\cdots,n\}$ and $x+(k+1)t>n$, then $k\ge 1$ and $x+\frac {k+1}{2}t>n/2\ge m$.
Thus  $|A\cap \{1,2,\cdots,m\}|\le (k+1)/2$ when $k$ is odd, or $|A\cap \{1,2,\cdots,m\}|\le k/2+1$ when $k$ is even.
Hence $|A\cap \{1,2,\cdots,m\}|/|A|\le 2/3$.
Since $t>3n/r_n$, there is $h\ge 3$, $0\le i_1<i_2<\cdots i_h<r_n$ such that
$b+(i_1+1)t>n\ge  b+i_1 t,  b+(i_2+1)t>2n\ge  b+i_2 t, \cdots, b+(i_h+1)t>hn\ge  b+i_h t$
and $b+r_nt\le (h+1)n$.
Let $i_0=-1$, $i_{h+1}=r_n$ and
$$A_j=\{b+(i_{j}+1)t \mod n, b+(i_{j}+2)t\mod n,\cdots, b+i_{j+1} t\mod n\}.$$
Then $|A_j\cap\{1,2,\cdots m\}|\le \frac 23|A_j|$ when $0\le j\le h-1$, and
$\jmath=|A_{h+1}\cap\{1,2\cdots,m\}|\le m/t+1$.
Thus
\beq
&&|\{b,b+t\mod n,\cdots, b+r_n t\mod n\}\cap \{1,2,\cdots, m\}|/r_n\\
&&\le
\frac{\frac 23 \sum_{j=0}^{h-1}|A_j|+\jmath}{r_n}
\le \frac {2(i_h+1+\jmath)}{3r_n}+\frac {\jmath}{3r_n}\le \frac 23+\frac m{3r_nt}+\frac 1{r_n}<\frac 34.\qed
\eeq

\textbf{Proof of Theorem 1.2}
(1)Firstly, by (\ref{s1-d2}) and (\ref{s3-l2-d2}),
\beqn
|e^{-\tilde I_{n,r}}-e^{-\mu_{n,r}}|\le |\tilde I_{n,r}- \mu_{n,r}|\le qnp^r r^2.\lbl{s3-pd1}
\eeqn
Next, let $r_n=[\frac {-2\ln n}{\ln p}+ \frac {\ln\ln n}{\ln p}+\frac {\ln 2-\ln q}{\ln p}-1]$
and
$R_n=[\frac {-3\ln n}{\ln p}]$.
Lemma \ref{s3-l3}, together with (\ref{s3-d3}), (\ref{s3-d4}), (\ref{s3-l2-d1}) and (\ref{s3-pd1}), yields that
\begin{eqnarray}
&&\max_{r_n\le r\le R_n}|P(W^{(n)}<r)-e^{-\mu_{n,r}}|\nn\\
&&\le \max_{r_n\le r\le R_n}|P(W^{(n)}<r)-e^{-\tilde I_{n,r}}|+\max_{r_n\le r\le R_n}|e^{-\tilde I_{n,r}}-e^{-\mu_{n,r}}|\nn\\
&&\le  \max_{r_n\le r\le R_n} (\tilde e^{(n,r)}+\frac {np^r}{qr}+ qnp^{r}r^2)\nn\\
&&= O(n^{-1}\ln^{7} n).\lbl{th2-pd100}
\end{eqnarray}
Furthermore, it's easy to check that
$e^{-\mu_{n,r_n}}=o(n^{-1})$
and
$1-e^{-\mu_{n,R_n}}=O(n^{-1})$.
Therefore(\ref{th2-d1}) holds by noting that
$P(W^{(n)}<r)$ and $e^{-\mu_{n,r}}$ are both increasing functions of $r$.

(2) Let $r=C\ln n+x$.
Then $\mu_{n,r}= {qp^x}/{2} $.
Hence (\ref{th2-d1}) implies that
\begin{eqnarray*}
&&\max_{x\ge D{\ln\ln n}+a,C\ln n+x\in \mathbb Z}|\exp( {qp^x}/{2} ) P(W^{(n)}<C\ln n +x)-1|\\
&\le & \exp(qp^{D \ln\ln n+a}/2) O(\frac {\ln^7 n}{n})=O(\frac{\ln^7 n}{n^{1-\frac {qp^a}{2}}}).
\end{eqnarray*}

(3) In view of  (\ref{th2-d1}), it holds that
\beq
P(W^{(n)}<C\ln n+\varepsilon D\ln\ln n)=e^{-O(\ln^\varepsilon n)}+O(\frac {\ln^7 n}{n})\rightarrow\left\{
                                                                                                                 \begin{array}{ll}
                                                                                                                   1, & \hbox{$\varepsilon<0$;} \\
                                                                                                                   0, & \hbox{$\varepsilon>0$.}
                                                                                                                 \end{array}
                                                                                                               \right.
\eeq
It follows (\ref{th2-d2}) immediately.

(4) In view of (\ref{d1}), to prove (\ref{d2}), it remains only to show that
$$\limsup_{n\to\infty}\frac {W^{(n)}}{C\ln n}\le \frac 32 \,\, a.s. $$
For any $\varepsilon>0$, by (\ref{s3-d3}),(\ref{s3-l2-d1}) and (\ref{s3-l2-d2}),
\beqn
P(W^{(n)}>(1+\varepsilon)C\ln n) \le O(n^{-2\varepsilon})+O(\frac{n^{-1-2\varepsilon}}{\ln n}).\lbl{s3-pd2}
\eeqn
Hence
$$\sum_{n=1}^\infty P(W^{(n)}>(1+\varepsilon)C\ln n) <\infty$$
whenever $\varepsilon>\frac 12$.
Therefore,
 $\limsup_{n\to\infty}\frac {W^{(n)}}{C\ln n}\le \frac 32  \,\,\,a.s$.
as desired.

Since $W^{(n)}\ge U^{(n)}$, by (\ref{th1-d4}),
to prove (\ref{th2-d4}), we need only to show that
\beqn
 \liminf_{n\to\infty} \frac {W^{(n)}-C\ln n}{\ln\ln n}\le D\,\,\, a.s. \lbl{th2-pd5}
\eeqn
Fix any $0<\varepsilon<1$.
Let $r_n=[C\ln n+\varepsilon D\ln\ln n]$,  $H_n=\{W^{(n)}<r_n\}$ and
$X_k=\sum_{n=k+1}^{2k} I_{H_n}$.
Then
\beqn
P(\bigcup_{n=k+1}^{2k}H_n)=P(X_k>0)\ge \frac{(EX_k)^2}{EX^2_k}=\frac {(\sum_{n=k+1}^{2k} P(H_n))^2}{\sum_{m,n=k+1}^{2n} P(H_mH_n)}. \lbl{th2-pd30}
\eeqn
Clearly, $\mu_{n,r_n}=O(\ln^\varepsilon n)=o(\ln n/4)$. Together with (\ref{th2-d1}), it implies that
\beqn
\sum_{n=k+1}^{2k} P(H_n)\ge \sum_{n=k+1}^{2k}e^{-\mu_{n,r_n}}-O(\ln^7k)\ge O(k^{3/4}). \lbl{th2-pd31}
\eeqn
Set
$E_n=\{(a,s):1\le a\le n, [\frac {3n}{C\ln n}]+1\le s\le [\frac n2],\gcd(n,s)<\frac {n}{r_n}\}$.
By (\ref{s3-d1}) and (\ref{s3-l1-d1}),
$$P((H_mH_n)^c)\ge P(\cup_{(a,s)\in E_m}\tilde A^{(m,r_m)}_{a,s}\cup_{(b,t)\in E_n}\tilde A^{(n,r_n)}_{b,t}).$$
Let $V$ be the graph with vertex set $V_{k}=\{(a,s,n):k+1\le n\le 2k, (a,s)\in E_n\}$ and edges defined by
$(a,s,m)\thicksim (b,t,n)$ if and only if $\tilde B^{(m,r_m)}_{a,s}\cap \tilde B^{(n,r_n)}_{b,t}\not=\emptyset$.
Then $V$ is a dependency graph of $\{I_{\tilde{A}^{(n,r_n)}_{a,s}}:(a,s,n)\in V_k\}$.
Let $J_{m}=\sum_{(a,s)\in E_m} P(\tilde{A}^{(m,r_m)}_{a,s})$. Then
\beqn
P(H_mH_n)\le e^{-J_m-J_n}+\tilde e^{(m,r_m)}+\tilde e^{(n,r_n)}+2 \tilde e^{(m,r_m,n,r_n)}, \lbl{th2-pd33}
\eeqn
where
$\tilde e^{(m,r_m,n,r_n)}=\sum_{(a,s,m)\thicksim (b,t,n)}\big(P(\tilde{A}^{(m,r_m)}_{a,s})P(\tilde{A}^{(n,r_n)}_{b,t})+P(\tilde{A}^{(m,r_m)}_{a,s}\tilde{A}^{(n,r_n)}_{b,t})\big)$.
It follows that
\beqn
\sum_{m,n=k+1}^{2k} P(H_mH_n)\le (\sum_{n=k+1}^{2k} e^{-J_n})^2+{\cal L}_k,  \lbl{th2-pd34}
\eeqn
where ${\cal L}_k=2k\sum_{n=k+1}^{2k}\tilde e^{(n,r_n)}+2\sum_{m,n=k+1}^{2k}\tilde e^{(m,r_m,n,r_n)} $.
One deduces from (\ref{s3-l2-d2})  that
$J_n\ge \mu_{n,r_n}(1-\frac {(r_n+1)^2}{2n}-\frac {6}{C\ln n})$.
Hence
\beqn
 \sum_{n=k+1}^{2k} e^{-J_n}\le e^{O(\ln^{\varepsilon-1}k)} \sum_{n=k+1}^{2k} e^{- \mu_{n,r_n}} \lbl{th2-pd35}
\eeqn
by noting that $\mu_{n,r_n}=O(\ln^\varepsilon n)$.
If we have showed that
${\cal L}_k\le O(k\ln^8 k)$, then
$\lim_{k\to\infty}P(\bigcup_{n=k+1}^{2k}H_n)=1$ and   hence (\ref{th2-pd5}) holds,
 in view of (\ref{th2-pd30})--(\ref{th2-pd35}).

Now we shall show that ${\cal L}_k\le O(k \ln^8 k)$. By Lemma \ref{s3-l3},
$$ k\sum_{n=k+1}^{2k}\tilde e^{(n,r_n)}\le O(k\ln^7k).$$
Fix any $(a,s,m)\in V_k$, define
$\gamma(a,s,m)$  to be the set of all triples $(b,t,n)\in V_k$ such that $|\tilde{B}^{(n,r_n)}_{b,t}\cap\tilde{B}^{(m,r_m)}_{a,s}|> {C\ln {(2k)}}/2$.
Suppose that $k$ is sufficiently large and  $(b,t,n)\in \gamma(a,s,m)$.
Let $i_l=\min\{j\ge 0:b+jt>ln\}$, $\jmath=\max\{l:i_l\le r_m\}$ and
$Z_l=\{b+it-ln:i_l\le i\le min(i_{l+1}-1,r_n)\}$. Then
$\tilde{B}^{(n,r_n)}_{b,t}$ is the disjoint union of $Z_i$ with $0\le i\le \jmath$.
Since  $[\frac {3n}{C\ln n}]+1\le t\le [\frac n2]$,
 $|Z_i|\le C\ln n/3+1$ for all $i\le\jmath$, and $|Z_i|\ge 2$ for $0<i<\jmath$.
Thus $|\tilde{B}^{(n,r_n)}_{b,t}|=r_n+1\ge 2(\jmath-1)$.
It follows that
$ \jmath< C\ln{ (2k)}/2-3$.
Since  $|\tilde{B}^{(n,r_n)}_{b,t}\cap\tilde{B}^{(m,r_m)}_{a,s}|> {C\ln {(2k)}}/2$ and $|Z_i|\le C\ln n/3+1$,
there are $l$ and $\imath$ such that $\imath\not=l$,  $|Z_l\cap\tilde{B}^{(m,r_m)}_{a,s}|\ge 2$
and $|Z_\imath\cap\tilde{B}^{(m,r_m)}_{a,s}|\ge 1$.
That is to say, there is
$0\le i,j,\ell\le r_n$ and $x,y,z\in \tilde{B}_{a,s,m}$ such that
$b+it-ln=x$, $b+jt-ln=y$ and $b+\ell t-\imath n=z$. It leads to $(j-i)t=y-x$
and $(\imath-l)n=(\ell-i)t+x-z$.
Therefore
\beq
|\gamma(a,s,m)|\le &&|\{(b,t,n): t=(y-x)/i,
 n=(\ell t+x-z)/j,  b+\imath t \mod n=x,\\
&&|i|,|j|,|\ell|,|\imath|\le C\ln 2k, i,j\not=0,x,y,z\in\tilde{B}_{a,s,m}\}|=O(\ln^7 k).
\eeq
Combining with (\ref{s3-l4-d1}) and (\ref{s3-l4-d2}),  one  then  deduces that
\beq
\sum_{m,n=k+1}^{2k}\tilde e^{(m,r_m,n,r_n)} &\le& \sum_{m,n=k+1}^{2k}\big(2p^{r_m+r_n-1}m^2nr_mr_n
+3p^{r_m+r_n-C\ln(2k)/2}m^2r^2_m r^3_n\big)\\
&&+\sum_{(a,s,m)\in V_k}p^{r_m}|\gamma(a,s,m)|\le O(k\ln^8 k)
\eeq
as desired.
It completes the proof of  (\ref{th2-d4}).

(5)
Now we come to the proof of (\ref{th2-d5}).
By (\ref{th2-d1}), for any $0<\varepsilon<1$,
\beq
P(W^{(2^n)}<C\ln 2^n+\varepsilon D\ln\ln 2^n)=e^{-O(n^\varepsilon)}
\eeq
and
\beq
P(W^{(2^n)}>C\ln 2^n-(1+\varepsilon) D\ln\ln 2^n)=O(n^{-(1+\varepsilon)}).
\eeq
Hence
$\liminf_{n\to\infty}\frac {W^{(2^n)}-C\ln 2^n}{\ln\ln 2^n}\ge 0$ and  $\limsup_{n\to\infty}\frac {W^{(2^n)}-C\ln 2^n}{\ln\ln 2^n}\le -D$ by
the Borel-Cantelli Lemma.
In view of (\ref{th2-d2}), it remains only to show that
\beqn
\limsup_{n\to\infty}\frac {W^{(2^n)}-C\ln 2^n}{\ln\ln 2^n}\ge -D.  \lbl{th2-pd41}
\eeqn
For any $0<\varepsilon<1$, let $r_n=[C\ln n-\varepsilon D\ln\ln n]$ and
$F_{n}=\{(a,s):1\le a\le n, 3n/r_n<s\le [n/2], \gcd(n,s)<n/r_n\}$.
To show (\ref{th2-pd41}), it suffices  to show that
$\lim_{n\to\infty}P(\bigcup_{n=k+1}^{2k} \{W^{(2^n)}\ge r_{2^n}\})=1$.
It need only to show that
$\lim_{n\to\infty}P(\bigcup_{(a,s,m)\in G_k}\tilde A^{(m,r_m)}_{a,s})=1$, where
$G_k=\{(a,s,2^n):k+1\le n\le 2k, (a,s)\in F_{2^n}\}$.
Let $X_k=\sum_{(a,s,m)\in G_k} I_{\tilde A^{(m,r_m)}_{a,s}}$.
Clearly $EX_k=O(k^{1-\varepsilon})$.
By Lemma \ref{s3-l3}, Lemma \ref{s3-l4}, (\ref{s3-l4-d1}) and (\ref{s3-l4-d2}),
\beq
DX_k &&\le \sum_{(a,s,m)\in G_k }P(\tilde A^{(m,r_m)}_{a,s})+2\sum_{n=k+1}^{2k} \tilde e^{(2^n,r_{2^n})}\\
&&+2\sum_{k+1\le m<n\le 2k} (p^{r_{2^m}+r_{2^n}-1}2^{2m}2^nr_{2^m}r_{2^n}+3p^{r_{2^m}+\frac 14 r_{2^n}}2^{2m}r^2_{2^m}r^3_{2^n})\\
&&= O(k^{1-\varepsilon}).
\eeq
Hence
$$P(\bigcup_{(a,s,m)\in G_k}\tilde A^{(m,r_m)}_{a,s})=P(X_k>0)\ge \frac{(EX_k)^2}
{DX_k+(EX_k)^2}\to 1.$$

As to (\ref{th2-d7}), by Proposition \ref{as-p3}, we need only to show
(\ref{as-d10}). It suffices to show that
\beqn
\lim_{n\to\infty} \frac 1n \sum_{k=1}^n \sum_{(a,s)\in F_{2^k}}I_{\tilde A^{(2^k,r_{2^k})}_{a,s}}=0 \,\, a.s. \lbl{th2-pd51}
\eeqn
 and
 \beqn
 \lim_{n\to\infty} \frac 1n \sum_{k=1}^n I_{H_k}=0 \,\, a.s.,\lbl{th2-pd52}
 \eeqn
 where $H_k=\{W^{(2^k)}\ge r_{2^k}\}\setminus (\bigcup_{(a,s)\in F_{2^k}}\tilde A^{(2^k,r_{2^k})}_{a,s})$.
Similar as the proof of (\ref{as-d6}), we may show (\ref{th2-pd51}).
Note that $P(H_k)=O(k^{-1-\varepsilon})$. By Borel-Cantelli Lemma,
$\lim_{k\to\infty}I_{H_k}=0 \,\,a.s$. It follows (\ref{th2-pd52}) and completes our proof.
\qed

\section{The asymptotic distributions when $p_n=o(1)$}
\setcounter{equation}{0}
\textbf{Proof of Theorem 1.3}
(i) Set $q_n=1-p_n$ and $r_n=[\frac {2\ln n}{-\ln p_n}+\frac {\ln\ln n}{\ln p_n}]$. Clearly,
 $p^{r_n+1}_n\le n^{-2}\ln n \le p^{r_n}_n$.
Similar as the proof of (\ref{th1-pd100}), we may show that there is a constant $c>0$ such that
\beqn
\max_{0\le  k\le 2}|P(U^{(n,p_n)}<r_n+k)-e^{-\lambda_{n,r_n+k,p_n}}|
\le c p^{-3}_nn^{-1}\ln^2 n. \lbl{th3-pd-1}
\eeqn
Since $\lim_{n\to\infty}\frac {2\ln n}{-\ln p_n}=\infty$, $\frac {\ln n}{-\ln p_n}\ge 10$ and hence
$p^{-1}_n\le n^{0.1}$ for sufficiently  large $n$. This, together with (\ref{th3-pd-1}) implies that
\beqn
\lim_{n\to\infty}\max_{0\le  k\le 2}|P(U^{(n,p_n)}<r_n+k)-e^{-\lambda_{n,r_n+k,p_n}}|=0. \lbl{th3-pd-2}
\eeqn
In addition,
\beqn
\lim_{n\to\infty} \lambda_{n,r_n,p_n}\ge \lim_{n\to\infty}\frac {q_n\ln n}{2(r_n-1)}=
 \lim_{n\to\infty}\frac {-\ln p_n}{4}=\infty \lbl{th3-pd-3}
\eeqn
 by noting that
$p_n\rightarrow 0$.
Similarly,
\beqn
  \lim_{n\to\infty} \lambda_{n,r_n+2,p_n}\le\lim_{n\to\infty} (\frac {p^2_n\ln^2 p_n}{8\ln n}-\frac {p_n\ln p_n}{4})=0. \lbl{th3-pd-4}
\eeqn
Therefore (\ref{th3-d5}) holds by (\ref{th3-pd-2})--(\ref{th3-pd-4}).
Use the same method, we may prove (\ref{th3-d6}) and  complete the proof of (i).

(ii)(iii):
Since $\lim_{n\to\infty} np_n=\infty$,
$\lim_{n\to\infty}P(\sum_{i=1}^n \xi^{(n)}_i\ge 2)=1$ and hence
\beqn
\lim_{n\to\infty}P(U^{(n,p_n)}\ge 2)=1.\lbl{th3-pd0}
\eeqn
Choose an $0<\varepsilon<0.1$ such that $[b-\varepsilon, b)\cup(b,b+\varepsilon)$ contains no integers. There is $N$ such that
 $\frac{b-\varepsilon}{2}< \ln n/(-\ln p_n)<\frac{b+\varepsilon}{2}$ for all $n>N$. It follows that for all $n>N$,
\beqn
p^{\frac{-b+\varepsilon}{2}}_n<n<p^{\frac{-b-\varepsilon}{2}}_n. \lbl{th3-pd1}
\eeqn
 By   (\ref{s3-d3}), (\ref{th3-pd1}) and Lemma \ref{s3-l2}, for $n>N$,
\beqn
P(W^{(n,p_n)}\ge [b]+1)\le \frac {p_n^{[b]+1-b/2-\varepsilon/2}}{q_n([b]+1])}+\frac {q_n p_n^{[b]+1-b-\varepsilon}}{2}\rightarrow 0.\lbl{th3-pd2}
\eeqn
Suppose that $r<b\le r+1$ where $r$ is a positive integer. By (\ref{th3-pd1}), $n^2p^r_n\ge p^{r-(b-\varepsilon)}_n\to\infty$.
Let $X_n=\sum_{(a,s)\in B^{(r)}_n}I_{A^{(r)}_{a,s}}$. By Lemma \ref{s2-l1} and
Lemma \ref{s2-l2}, $EX_n=O(n^2p^r_n)$ and $DX_n\le EX_n+O(n^3p^{2r-1}_n)=o\big((EX_n)^2\big)$. Consequently,
\beqn
P(U^{(n,p_n)}\ge r)=P(X_n\ge 0)\ge \frac {(EX_n)^2}{DX_n+(EX_n)^2}\to 1.\lbl{th3-pd6}
\eeqn
Similar as the proof of (\ref{th1-pd100}) and (\ref{th2-pd100}), by using  (\ref{th3-pd1}),we may show that
\beqn
\lim_{n\to\infty}|P(U^{(n,p_n)}<b)-e^{-\lambda_{n,b,p_n}}|=\lim_{n\to\infty}|P(W^{(n,p_n)}<b)-e^{-\mu_{n,b,p_n}}|=0 \lbl{th3-pd3}
\eeqn
when  $b\ge 3$ and $b$ is an integer.
Furthermore, if $u=\lim_{n\to\infty} n^2p^b_n\le \infty$ exists, then
\beqn
\lim_{n\to\infty}P(U^{((n,p_n)}<b)=\lim_{n\to\infty}e^{-\lambda_{n,b,p_n}}=e^{-\frac {u}{2(b-1)}}\lbl{th3-pd4}
\eeqn
and
\beqn
\lim_{n\to\infty}P(W^{(n,p_n)}<b)=\lim_{n\to\infty}e^{-\mu_{n,b,p_n}}=e^{-\frac {u}{2}}.\lbl{th3-pd5}
\eeqn
Thus our result holds by (\ref{th3-pd0}),(\ref{th3-pd2}),(\ref{th3-pd6}),(\ref{th3-pd4}), (\ref{th3-pd5})
and by noting that $W^{(n,p_n)}\ge U^{(n,p_n)}$.
\qed

\setcounter{equation}{0}

\end{document}